# Data-Driven Online Optimization for Enhancing Power System Oscillation Damping

Zhihao Chen, *Student Member, IEEE,* Hanchen Xu, *Student Member, IEEE,* Junbo Zhang, *Senior Member, IEEE,* and Lin Guan, *Member, IEEE*

***Abstract*—This paper reports an initial work on power system oscillation damping improvement using a data-driven online optimization method. An online oscillation damping optimization model is proposed and formulated in a form solvable by the data-driven method. Key issues in the online optimization procedures, including the damping sensitivity identification method, its compatibility with the dispatch plans, as well as other practical issues in real large-scale system are discussed. Simulation results based on a 2-area 4-machine system, an 8-machine 36-bus system and NETS-NYPS 68-bus system verify the feasibility and efficiency of the proposed method. The results also show the capability of the proposed method to bridge the gap between online data analysis and complex optimization for power system dynamics.**

***Index Terms*—data-driven method, damping improvement, inter-area mode oscillation, online optimization, sensitivity identification.**

## I. Introduction

MODERN power systems are undergoing a series of evolutions including the popularization of power electronic devices, the deepening integration of massive renewable energy sources (RESs), the rapid development of electricity markets, and the extensive application of various heterogeneous controllers [1-3]. Under such circumstances, high uncertainty and complexity become new features of modern power systems, which increasingly complicates system dynamics.

One of the most complicate issues frequently encountered in modern power systems is the oscillations. In general, power system oscillations are caused by insufficient damping, forced oscillations, unmatched controller settings, unmatched network parameters, or nonlinear phenomenon [4-10]. Although most oscillations can be eliminated by controllers, it is still desirable that certain oscillations that are highly related to the system operating points, such as inter-area low frequency oscillations, can be resolved in an operation manner, rather than employing more sophisticated wide area controllers, which may cause even more control stability issues since the control systems may become more and more complex.

This work was supported in part by the National Natural Science Foundation of China (51607071), and the Natural Science Foundation of Guangdong Province, China (2018B030306041).

Zhihao Chen, Junbo Zhang (*Corresponding Author*) and Lin Guan are with School of Electrical Power, South China University of Technology, Guangzhou, Guangdong, P. R. China. E-mail: epjbzhang@scut.edu.cn.

Hanchen Xu is with Department of Electrical and Computer Engineering, University of Illinois at Urbana-Champaign, Urbana, IL 61801, USA.

To prevent system from inter-area low frequency oscillations, the system should be operated at a secure operating point, at which the minimum damping ratio of the critical inter-area oscillation mode is higher than 3-4%. However, in some cases, to mitigate the impacts from uncertainty caused by the RESs, the system may operate at a poorly damped point, at which the minimum damping ratio is lower than 3-4%, and then re-dispatch measures must be taken to restore the system to security again, which raises the problem of damping improvement using re-dispatch strategies [11-12].

Conventionally, damping improvement using re-dispatch strategy is designed offline with model-based optimization method [11-12], and it can hardly adapt to various operation conditions in the uncertain modern power systems. In recent years, with the development of signal processing and data science, online oscillation damping monitoring becomes possible [13-14], making it feasible to establish the connection between the oscillation damping and various system operation conditions, and consequently makes the data-driven based damping improvement feasible.

Some pioneers in this area have tried using damping sensitivity as an intermediate variable for system re-dispatch, but they ignored the fact that the damping sensitivity varies significantly along the re-dispatch process; therefore, these methods are generally inappropriate for real power systems [15]. Some other researchers tried artificial neural networks, but initial results were limited to a fixed system operation condition, making the algorithm impractical in real power systems [16]. In fact, while improving system damping with re-dispatch may sound attempting, developing a practical data-driven method for it is quite challenging due to 1) the time-varying nature of the optimization problem, 2) the accuracy problem of parameters used in the optimization, and 3) the poor data quality and the high cost for data communication [17-20]. Indeed, all these issues may lead to various practical concerns that need extensive investigation. In addition, machine learning algorithms for a large-scale system operation may also encounter issues such as incomplete data, underdetermined model, and ill-conditioned problems, which highlight the need of feature reduction using problem specific decomposition. In this paper, an initial work on power system oscillation damping improvement using data-driven based online optimization method is reported, which can also be regarded as a bridge between the situation awareness and practical power system operations and an attempt to justify the capability of online real-time data analysis in complicated power system dynamic optimizations. The proposed method has the following merits:

- purely data-driven without relying on system parameters and offline models;
- complementary to the damping controllers to solve low damping problems caused by fast load/renewable energy changing;
- online adaptable to various system operating conditions; and
- feasible to be implemented in large-scale systems using problem specific feature reduction techniques.

The remainder of this paper is organized as follows. Section II proposes a data-driven online optimization method for power system oscillation damping improvement. Section III discusses key issues in the optimization procedures, including the damping sensitivity identification and its compatibility to the dispatch plans, as well as practical issues for online applications. Section IV provides three case studies on a 2-area 4-machine system, and an 8-generator 36-node system as well as the NETS-NYPS 68-bus system, respectively. Section V concludes the paper.

## II. Data-Driven Online Damping Optimization

In this section, an initial model of damping optimization is first developed, and a modified model is presented, which can be solved by a data-driven approach online.

### A. Damping Optimization Problem

Power system operation conditions can be characterized by system states and parameters [4]. The former include power flows, device connection and operation states, while the latter include network parameters, device parameters, and etc. In general, when the system operates in a steady state, these parameters vary smoothly, and the system operation condition is mainly determined by the power flow, which is further determined by the active and reactive power injections at each bus in the system. Therefore, if we measure the power injections at each bus as $\boldsymbol{x}$, the system state matrix $\boldsymbol{A}$ of the system linearized model becomes a function of $\boldsymbol{x}$. Then, the damping ratios of the eigenvalues of $\boldsymbol{A}$, denoted by $\boldsymbol{\zeta}$, can be calculated as a function of $\boldsymbol{x}$ as follows:

$$\boldsymbol{\zeta} = \boldsymbol{f}(\boldsymbol{x})\,. \tag{1}$$

In this case, the online damping ratio optimization can be summarized as the following problem: given $\boldsymbol{x}[t]$ at time $t$, if the minimum damping ratio of the system, $\zeta[t]$, is lower than 3%, the objective is to determine the optimal changes of $\boldsymbol{x}$ at time $t$+1, $\Delta\boldsymbol{x}[t+1]$, so as to maximize the minimum damping ratio of $\zeta[t+1]$ determined by $\boldsymbol{x}[t+1]$:

$$\begin{aligned} &\max_{\Delta x}\ \min_{j}\{\boldsymbol{f}_j(\boldsymbol{x}[t]+\Delta\boldsymbol{x}[t+1])\} \\ &s.t.\ \ \Delta\underline{\boldsymbol{x}}[t+1]\leqslant\Delta\boldsymbol{x}[t+1]\leqslant\Delta\overline{\boldsymbol{x}}[t+1] \\ &\qquad \underline{\boldsymbol{x}}[t+1]\leqslant\boldsymbol{x}[t]+\Delta\boldsymbol{x}[t+1]\leqslant\overline{\boldsymbol{x}}[t+1] \end{aligned} \tag{2}$$

where $\Delta\underline{\boldsymbol{x}}[t+1]$ and $\Delta\overline{\boldsymbol{x}}[t+1]$ are the lower and upper bounds of $\Delta\boldsymbol{x}[t+1]$, respectively; $\underline{\boldsymbol{x}}[t+1]$ and $\overline{\boldsymbol{x}}[t+1]$ are the lower and upper bounds of $\boldsymbol{x}[t+1]$, respectively.

Note that, the lower and upper bounds of $\boldsymbol{x}$ and $\Delta\boldsymbol{x}$ at time $t$ are determined by the boundaries of generation, the generation ramp rates, the N-1 stability constraints, and etc. In this work, due to the space limitation, we assume the boundaries are already available and focus on the framework of the data-driven online damping optimization.

### B. Online Optimization Algorithm

The above problem can be solved when the system model $\boldsymbol{f}$ is available at each time $t$ and the computation time from $t$ to $t$+1 is sufficient. However, these two requirements are hardly satisfied because 1) the multi-machine bulk power system are integrated with intermittent and uncertain RESs and loads, and 2) the system has complex control networks which are difficult to model. Moreover, since the system is time-varying, it is impossible to obtain the real optimal point because the regulation takes time. Therefore, we modify the optimal problem so as to get a model that is applicable in real-time and gives near optimal results. Given $\boldsymbol{x}[t]$ at time $t$, approximating $\boldsymbol{f}$ with the following affine function:

$$\boldsymbol{f}_j(\boldsymbol{x}[t]+\Delta\boldsymbol{x}[t+1]) = \zeta[t]+\boldsymbol{\Psi}_j[t]\Delta\boldsymbol{x}[t+1]\,, \tag{3}$$

where $\boldsymbol{\Psi}[t]$ is the sensitivity matrix. Then, (2) can be formulated as:

$$\begin{aligned} &\max_{\Delta x}\ \zeta[t]+\boldsymbol{\Psi}_j[t]\Delta\boldsymbol{x}[t+1] \\ &s.t.\ \ \Delta\underline{\boldsymbol{x}}[t+1]\leqslant\Delta\boldsymbol{x}[t+1]\leqslant\Delta\overline{\boldsymbol{x}}[t+1] \end{aligned}\,. \tag{4}$$

Introducing an auxiliary variable $\delta$, the above optimal function can be further formulated as:

$$\begin{aligned} &\max_{\Delta x,\,\delta}\quad \delta \\ &s.t.\quad \boldsymbol{\Psi}_j[t]\Delta\boldsymbol{x}[t+1]\geqslant\delta \\ &\qquad \Delta\underline{\boldsymbol{x}}[t+1]\leqslant\Delta\boldsymbol{x}[t+1]\leqslant\Delta\overline{\boldsymbol{x}}[t+1] \end{aligned} \tag{5}$$

Problem (5) is a linear program (LP) and can be solved efficiently using a variety of algorithms [21].

### C. The Data-Driven Approach

Optimization problem (5) is solvable online if parameters $\boldsymbol{\Psi}[t]$, $\Delta\underline{\boldsymbol{x}}[t+1]$, and $\Delta\overline{\boldsymbol{x}}[t+1]$ are available. For power systems, $\Delta\underline{\boldsymbol{x}}[t+1]$ and $\Delta\overline{\boldsymbol{x}}[t+1]$ are mainly determined by the ramp rate and capability of generators, the dispatch intervals, and the original dispatch schedule [22]. While for $\boldsymbol{\Psi}[t]$, which is a nonlinear and time-varying vector, it can be estimated using the operation data provided by the wide-area measurement systems (WAMS) and the supervisory control and data acquisition system (SCADA). Therefore, we have the data-driven approach in three parallel threads, which is illustrated via the flow chart shown in Fig. 1.

Parallel thread 1: measure system states with WAMS and SCADA;

Parallel thread 2: estimate minimum damping of the system using a variety of methods such as the ones in [13] and [14];

Parallel thread 3: if the minimum damping ratio of the system, $\zeta[t]$, is lower than 3% at time $t$, estimate damping sensitivity $\boldsymbol{\Psi}[t]$ using the method to be discussed in Section III, and then solve (5) to re-dispatch the system at $t$+1 with the optimal result $\Delta\boldsymbol{x}[t+1]$.

Note that, among these three threads, the sensitivity matrix identification is just one part of the approach frame, but it is the very key to a successful optimization at a power system operating point.

In addition, it is possible that the minimum eigenvalue may change during the operation. However, the proposed method attempts to improve the damping ratio of the small signal stability problem but not to get trace of the same eigenvalue. In this case, the proposed method will adjust the minimum eigenvalue as long as the minimal one reach the requirement.

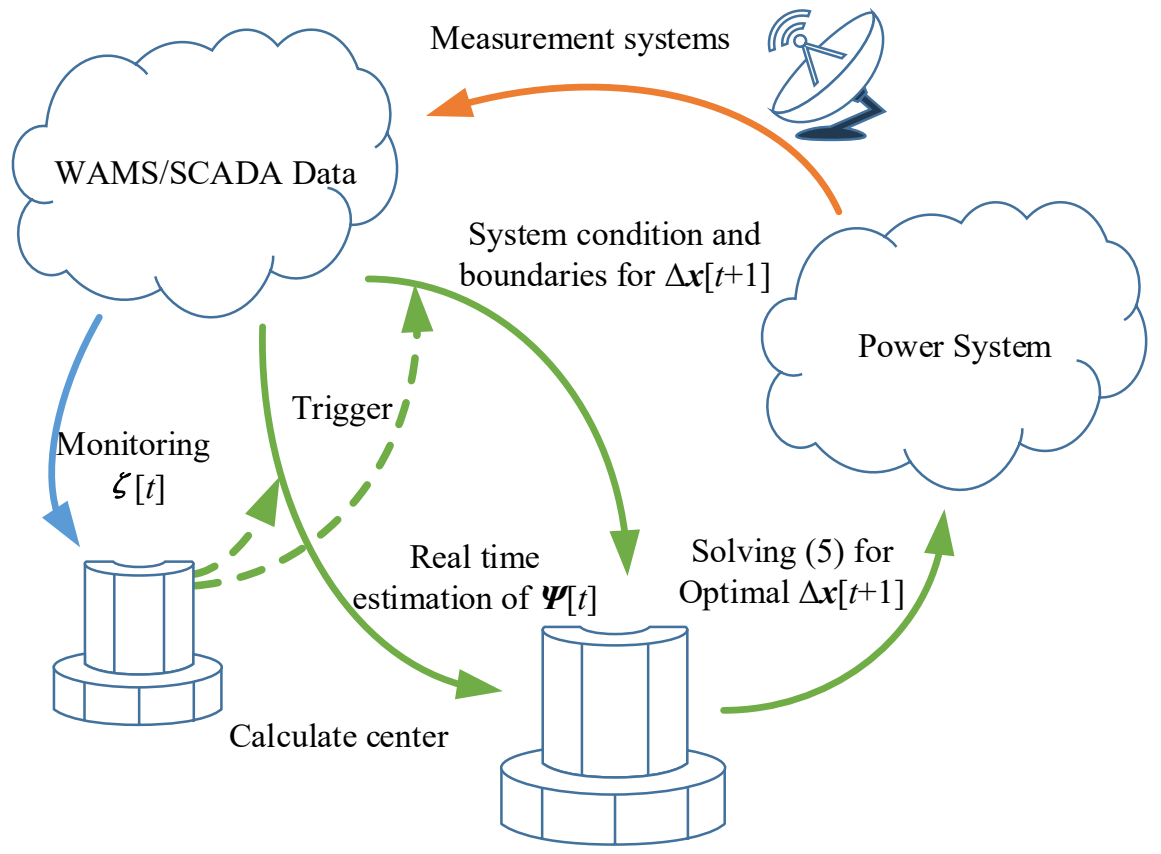

Fig. 1. Flow chart of the data-driven based online damping optimization.

### D. Discussions

To implement the proposed method in real power systems, there are some practically issues to be considered:

First, the matrix $\boldsymbol{A}$ of a power system is not only influenced by injections but also by grid topological modifications. The grid topological modifications will change the impedance in the direct path between the main groups of oscillating machines of an inter-area mode and therefore modify the damping of the mode. However, in large power systems, one or two grid topological modifications would not bring large changes to the damping ratio and the damping sensitivities of the dominant interarea mode since one grid topological modification is too limited. In this case, the proposed method can be reset according to the new operating data, and then it can work after the transient.

Second, the low damping problem in a real power system may cause by a fault or fast load/renewable energy change. In the former case, damping controllers such as power system stabilizers (PSSs) are equipped to increase the damping. While in the second case, generator re-dispatch should be performed, and the proposed approach could play a complementary role to the damping control.

Third, the sensitivity matrix identification is one part of the approach frame, but it is the very key to a successful optimization, which guarantees the damping ratio of the minimum eigenvalue will be improved at time $t$ compared to the case without optimization, and we are going to discuss this issue in the following section. However, it does not guarantee the minimum eigenvalue at time $t$ will always be improved compared to the one at time $t$-1 due to the planned dispatch and load/RESs change.

## III. Key Issues in Optimization Procedures

In this part, we provide solutions to the key issues in the online optimization procedures, including the damping sensitivity identification, the compatibility with the dispatch plans, and its potential to be implemented in large-scale power systems. Some other practical issues are also discussed.

### A. Damping Sensitivity Identification

Given $\boldsymbol{x}[t]$ at time $t$, the linearized form of (1) for the minimum damping mode is:

$$\Delta\zeta[t] = \boldsymbol{\Psi}_j[t]\Delta\boldsymbol{x}[t] + \varepsilon_j[t], \tag{6}$$

where $\boldsymbol{\Psi}_j[t]$ is a constant vector and $\varepsilon_j[t]$ is the error term which is independent and identically Gaussian distributed [23-24].

Suppose the system is stable and $\boldsymbol{x}[t]$ is an equilibrium operating point, (6) may hold for the same damping sensitivity $\boldsymbol{\Psi}_j$ at $t$ = 1, 2, ..., $N$, which leads to the following matrix form:

$$\Delta\boldsymbol{\zeta} = \boldsymbol{\Psi}_j\Delta\boldsymbol{X} + \boldsymbol{\varepsilon}_j, \tag{7}$$

where

$$\begin{aligned}\Delta\boldsymbol{\zeta} &= [\Delta\boldsymbol{\zeta}[1], \Delta\boldsymbol{\zeta}[2], ..., \Delta\boldsymbol{\zeta}[N]]^{\mathrm{T}},\\ \Delta\boldsymbol{X} &= [\Delta\boldsymbol{x}[1], \Delta\boldsymbol{x}[2], ..., \Delta\boldsymbol{x}[N]]^{\mathrm{T}},\\ \Delta\boldsymbol{\varepsilon}_j &= [\Delta\varepsilon_j[1], \Delta\varepsilon_j[2], ..., \Delta\varepsilon_j[N]]^{\mathrm{T}}.\end{aligned} \tag{8}$$

Since $\Delta\boldsymbol{\zeta}$ can be estimated with many ambient signal analysis methods, $\Delta\boldsymbol{X}$ is directly measurable, $\boldsymbol{\Psi}_j$ is solvable with regression methods.

Generally, to fit the nonlinear and time-varying power system property, as well as to overcome collinearity and noise impacts of the system data, a noise assisted ensemble regression is applied. We can estimate $\boldsymbol{\Psi}_j$ with the cost function in (9), the estimator in (10).

$$J = \left\|\Delta\hat{\zeta} - \Delta\zeta\right\|_W^2 + k\left\|\hat{\boldsymbol{\Psi}}_j\right\|^2 \tag{9}$$

$$\hat{\boldsymbol{\Psi}}_j = (\Delta\boldsymbol{X}^T\boldsymbol{W}\Delta\boldsymbol{X} + k\boldsymbol{I})^{-1}\Delta\boldsymbol{X}^T\boldsymbol{W}\Delta\boldsymbol{\zeta} \tag{10}$$

where $\boldsymbol{W}$ = diag$\{w_i\}$ is a diagonal weight matrix.

Some parameters appeared in the algorithm can be tuned with the same methods in [25] and [26]. To focus on the most important issues and also due to the space limitation, we do not provide more discussions here. We refer the interested readers to the above references for more details.

### B. Re-dispatch with Compatibility to the Planned Dispatch

In order to be compatible with other planned dispatch applications, we have the following considerations:

1) Given $\boldsymbol{x}[t]$ at time $t$, the dispatched value at time $t$+1 should be formulated as $\boldsymbol{x}[t+1] = \boldsymbol{x}[t] + \Delta\boldsymbol{x}[t+1]_{\mathrm{R}} + \Delta\boldsymbol{x}[t+1]_{\mathrm{O}}$, where $\Delta\boldsymbol{x}[t+1]_{\mathrm{R}}$ is the solution of (5) and $\Delta\boldsymbol{x}[t+1]_{\mathrm{O}}$ is the dis-patch value for other planned dispatch such as day ahead planning and economic dispatch [23].

2) Constraint $\Delta\underline{\boldsymbol{x}}[t+1] \leqslant \Delta\boldsymbol{x}[t+1] \leqslant \Delta\overline{\boldsymbol{x}}[t+1]$ in (5) is modified as $\Delta\underline{\boldsymbol{x}}[t+1] \leqslant \Delta\boldsymbol{x}[t+1]_{\mathrm{R}} + \Delta\boldsymbol{x}[t+1]_{\mathrm{O}} \leqslant \Delta\overline{\boldsymbol{x}}[t+1]$, where $\Delta\boldsymbol{x}[t+1]_{\mathrm{O}}$ is considered as a known constant. $\Delta\underline{\boldsymbol{x}}[t+1] = \max\{-\Delta\boldsymbol{x}_{\mathrm{Ramp}}, \underline{\boldsymbol{x}}[t+1]-\boldsymbol{x}[t]\}$ and $\Delta\overline{\boldsymbol{x}}[t+1] = \min\{\Delta\boldsymbol{x}_{\mathrm{Ramp}}, \overline{\boldsymbol{x}}[t+1]-\boldsymbol{x}[t]\}$, where $\Delta\boldsymbol{x}_{\mathrm{Ramp}}$ is the limitation by the ramp rate.

3) Once a generator hit the limit, its re-dispatch value $\Delta\boldsymbol{x}[t+1]_{\mathrm{R}}$ becomes zero while its planned dispatch value $\Delta\boldsymbol{x}[t+1]_{\mathrm{O}}$ should be distributed to other generators.

The above considerations are integrated into thread 3 in Section II when damping optimization is online processed.

### C. Potential to Large-Scale Power System Implementation

Therefore, the implementation of the proposed method in large-scale power systems is an important concern that needs to be tackled carefully. The major challenge lays on the explosion of the feature space of the input-output regression model in (7), where variables $\Delta\boldsymbol{X}$ and $\Delta\boldsymbol{\zeta}$ could be thousands of dimensions. In this case, directly solving (7) is impossible. However, we can focus on reducing the feature space of the problem using our domain knowledge.

Basically, there are two ways achieving the above goal, one is from the side of the input variable space, and the other is from the side of the output variable space. Fortunately, the latter is quite easy to be accomplished because oscillation modes with poor damping are quite limited and obvious in many large power systems [27]. For example, the dominant mode of the ENTSO-E system is the one dominates the oscillation between generators in Turkey, Spanish and Portugal [28-29]. Therefore, once the dominant mode is selected, the target is to narrow the input feature space, which is also tricky enough because: 1) a mode in a system is in general dominated by a part of generators, say less than 1/3 of the total population in a common sense, which significantly reduces the input feature space; 2) the automatic generation control is implemented to power stations with several to tens of generators, especially for the hydro-power stations, which means the input features could be further divided by a factor of several to tens; 3) dispatching center is often layer-structured, meaning commands from the top layer are generally delivered to areas/districts, therefore small generators such as wind farms and PV stations are non-dispatchable and treated as irrelevant inputs. Given these tricks, any large power system with thousands of generators could be feature reduced using this problem specific decomposition by a factor of tens to a hundred, making the input space limited in a solvable scale, and consequently (7) is solvable.

It should be noticed that, feature reduction is quite different from model reduction or dynamic model equivalent. The former removes irrelevant generators and clusters dispatchable generators in one station only in the algorithm part for the data processing; while the latter aggregates generators and change the simulation model.

### D. More Discussions for the Online Applications

A typical illustration of thread 3 is shown in Fig. 2. $T_1$ is the time interval for re-dispatch and $T_2$ is the actual time it takes for the re-dispatch. Usually, $T_1$ is a constant while $T_2$ varies with $\Delta\boldsymbol{x}_{\mathrm{R}} + \Delta\boldsymbol{x}_{\mathrm{O}}$ at each time $t$.

Generally speaking, the damping sensitivity is nonlinear and varies along with the system operation condition. Therefore, at each time $t$, to estimate the sensitivity at $\boldsymbol{x}[t]$, the generation and load curves should not change much. This condition is usually accomplished by the load, because $T_1$ is usually less than 15 minutes, and the load change is usually not significant during that time interval. However, for generations with the re-dispatch value $\Delta\boldsymbol{x}_{\mathrm{R}} + \Delta\boldsymbol{x}_{\mathrm{O}}$, sometimes it may be difficult to meet this condition. Therefore, the measured data for the generation curves at time $t$ should be within the time interval from $t$-1+$T_2[t]$ to $t$.

In the previous work in [30], it has been pointed out the sensitivity identification has a prerequisite called the sufficient effective data condition that the power system sensitivity around an equilibrium operating point is identifiable online if the effective data around the point are sufficient and anisotropically dispersed, which guarantees the local property of the nonlinear model around the point can be approximated by a local affine plane. Therefore, to satisfy this prerequisite within a time interval from $t$-1+$T_2[t]$ to $t$, the sampling rate for the samples in (7) should be enhanced to its maximum, and the regression variables in $\Delta\boldsymbol{X}$ should be limited to these with high impact to the inter-area mode.

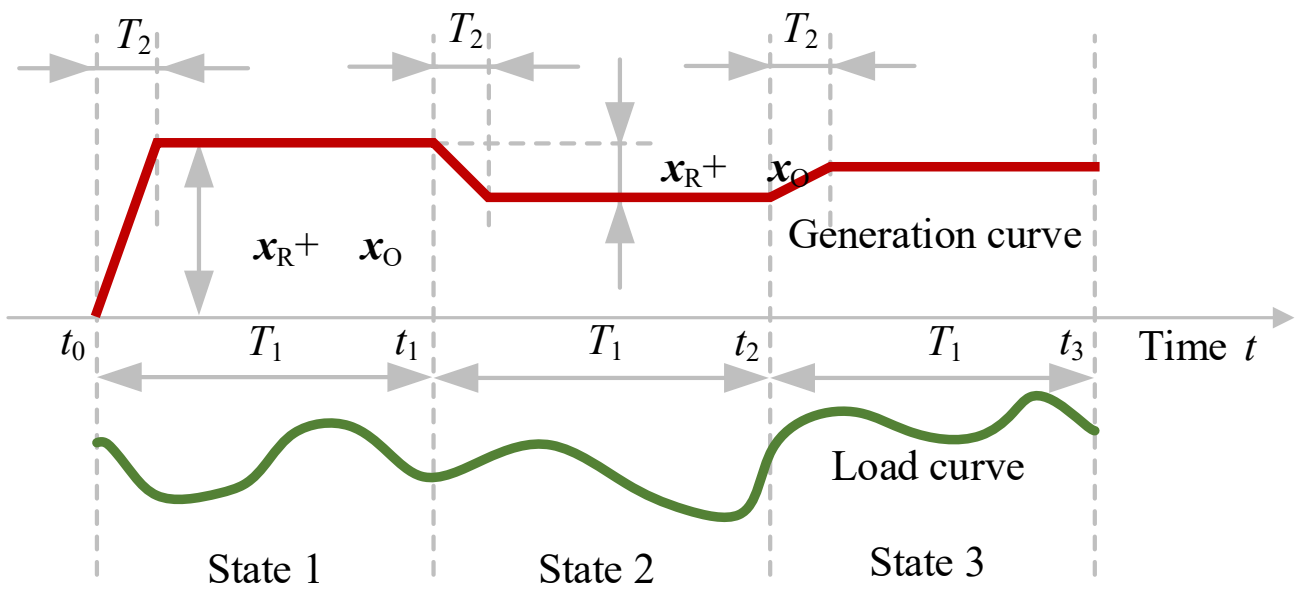

Fig. 2. Illustration for thread 3.

## IV. Case Study

In this section, the proposed method is tested in three benchmark systems for small signal stability analysis. Due to the space limitation, we provide detailed simulation information on the case with the smaller system, and then apply the same methodology to the larger systems to test its potential applicability to real power systems. The simulation is carried out in Matlab with the third-party software Power System Analysis Toolbox [31].

### A. Two-Area Four-Machine System

The two-area four-machine system in [4] is used, the diagram of which is illustrated in Fig. 3. The system has one inter-area mode with Generators 1 and 2 in Area 1 swinging against Generators 3 and 4 in Area 2, and two local modes in Area 1 and Area 2. In the initial system condition, no excitation system is equipped on the generators, the inter-area mode is unstable, and the local modes are stable. To achieve a stable condition, each generator is equipped with a casual type excitation as in [4], with a suitable gain to create a poorly damped inter-area mode, while the local modes are very stable. In the following parts, we focus on improving damping of the inter-area mode with the proposed method.

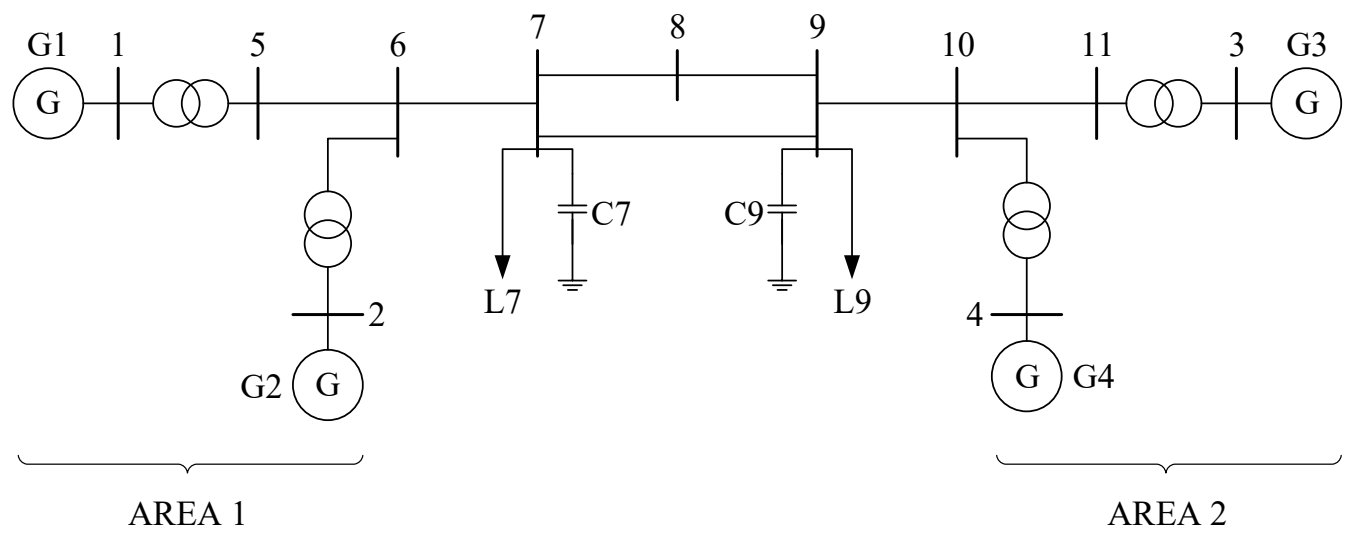

Fig. 3. The two-area four-machine system.

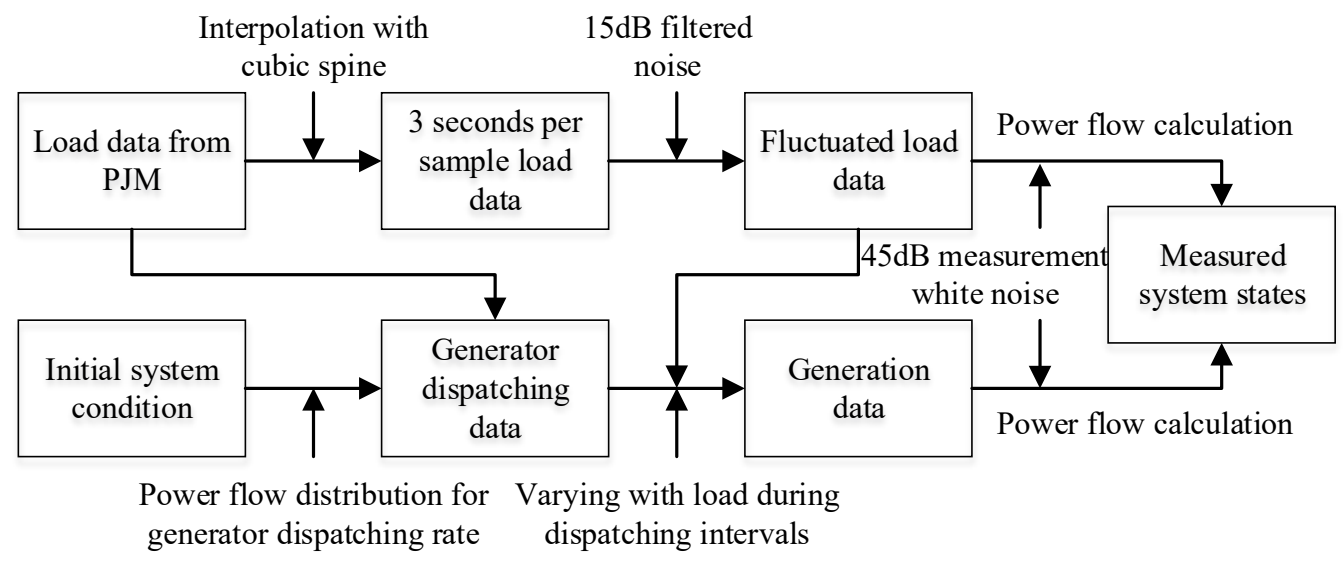

Fig. 4. Data generation procedure in P. 1.

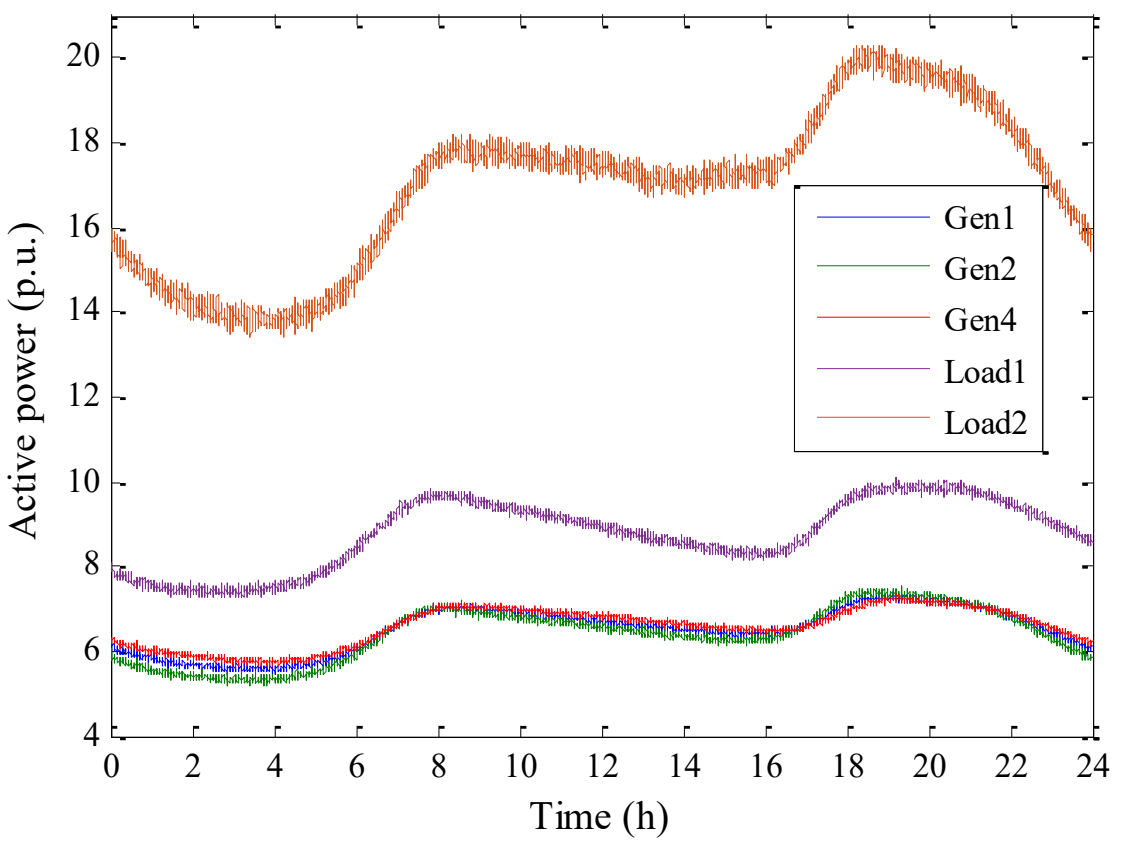

Fig. 5. Active power flow of the system before optimization.

TABLE I

THE ESTIMATED $\boldsymbol{\Psi}$[1] AND THE ONE WITH MODEL BASED METHOD

| Method | Gen1 | Gen2 | Gen4 | Load1 | Load2 |
|---|---|---|---|---|---|
| NAER | 0.0003 | 0.0033 | 0.0035 | -0.0058 | -0.0059 |
| Model based method | 0.0010 | 0.0032 | 0.0036 | -0.0057 | -0.0063 |

The three parallel threads in Section II Part C are employed.

1) Parallel thread 1:

Power flow calculation is carried out under different power system operation conditions. All generators except Generator 3, which is the slack bus, are considered as P-V buses, the loads are P-Q buses. The data generation is shown in Fig. 4.

To imitate real power system operations, the apparent power of the loads is varied along with time. One day of load data in 2013 is taken from the Pennsylvania-New Jersey-Maryland (PJM) Interconnection [32]. The data have two load peaks at about 8 a.m. and 7 p.m. To simulate load fluctuations, we interpolated the data with a cubic spline algorithm with sampling interval 3 seconds, and then added filtered white noise to the interpolated data with a signal-to-noise ratio (SNR) of 15dB. Active power of the generators is dispatched every 15 minutes according to the varied load data, with the same dispatch rate as in the initial condition. During the dispatch intervals, the generation varies with load fluctuations and Generator 3 is employed for power balance. The generation and load curves are then used for power flow calculation in PSAT for system states under different operating points. After that, about 45dB measurement white noise is added to the states. The final generation and load curves ($\boldsymbol{x}[t]$) are shown in Fig. 5.

2) Parallel thread 2:

Damping of the inter-area mode varies along with the power flow curves. It could be estimated using ambient signal based methods in [13-14]; to simplify the data generation in this study and to focus on the validation of the framework, it is calculated using the model-based method in PSAT. Then about 45dB white noise is added to the initial damping curve to manipulate the bias in the ambient signal based estimation.

The final damping curve ($\zeta[t]$) is shown in Fig. 6. It shows damping of the inter-area mode is very poor at the beginning operation point. It increases with the decreasing loads before 4 o'clock in the morning, and then decreases significantly to near zero with the first load peak at about 8 a.m. After that, it increases slowly when the loads are reduced smoothly. When the second load peak comes at about 5 p.m., the damping ratio decreases sharply and even reduces to negative values during the high load operations, indicating the system becomes unstable. After that, the system becomes stable again when the load peak recedes. Clearly, the proposed method should improve the damping to 3% before the first load peak and guarantee a stable system during the second load peak.

3) Parallel thread 3:

The initial damping of the inter-area mode is 1.75%, which is lower than 3%, and therefore, according to the trigger of parallel thread 3, the online optimization is carried out at the subsequent point of time.

In the first constraint in (5), the damping sensitivity $\boldsymbol{\Psi}[t]$ is estimated using the noise assisted ensemble regression in Section III Part A. Due to the manipulation on load fluctuations and measurement noise, the data $\boldsymbol{x}[t]$ does not have serious collinearity issue, the norm-2 penetration factor $k$ in (9) is set to 0, the number of ensemble is 100. With some matrix manipulations as in [25] and [26], $\boldsymbol{\Psi}[t]$ is estimated. Table I shows $\boldsymbol{\Psi}[1]$ for an example, where the estimated sensitivity is similar to the one calculated by the model based perturbation method. For $t$ equals to other time stamps, the sensitivity estimation has similar performance, which validates the effectiveness of the noise assisted ensemble regression.

In the second constraint in (5), $\Delta\underline{\boldsymbol{x}}[t+1] = \max\{-\Delta\boldsymbol{x}_{\text{Ramp}},$

$\underline{\boldsymbol{x}}[t+1]-\boldsymbol{x}[t]\ \}$, $\Delta\overline{\boldsymbol{x}}[t+1]\ =\ \min\{\ \Delta\boldsymbol{x}_{\mathrm{Ramp}}\ ,\ \overline{\boldsymbol{x}}[t+1]-\boldsymbol{x}[t]\ \}$, and the value of $\Delta\boldsymbol{x}[t+1]_{\mathrm{O}}$ is calculated as the dispatch increment of the generation curves in Fig. 5. $\overline{\boldsymbol{x}}[t+1]$ is the generator capacity, $\underline{\boldsymbol{x}}[t+1]$ is set to 20% of the generator capacity, and $\Delta\boldsymbol{x}_{\mathrm{Ramp}}$ is 5% of the generator capacity.

With all parameters known, the optimization (5) is then solved sequentially. The improved damping curve for the inter-area mode is shown in Fig. 6, and the active power after the re-dispatch is shown in Fig. 7. These results show that, two rounds of sequential optimization are carried out during the day. The first round is before 1:00 am, with mainly 4 steps of sequential optimization (see details in the first circle of Fig. 7); and the second round is during the first load peak, with mainly 7 steps of sequential optimization (see details in the second circle of Fig. 7). With the re-dispatch strategy, damping of the inter-area mode is improved to 3% before 1:00am and the system is stable during the second load peak. However, there is no re-dispatch during the second load peak even when the damping is lower than 3%. To find out what happens during that time, more information about the optimization is investigated.

The variations of damping sensitivity during the two rounds of sequential optimization are shown in Fig. 8. The left plot shows the variation of damping sensitivity during the first round, while the right plot shows those during the second round. In the first round, the sensitivity results suggest reducing active power of Generator 1 while increasing active power of Generator 4, which is consistent with the optimization based re-dispatch in Fig. 7; after 4 steps of sequential optimization, the damping is over 3% (Fig. 6) and then the re-dispatch stops temporarily. In the second round, the sensitivity results still suggest reducing active power of Generator 1 while increase that of Generator 4 (Fig. 8); however, Generator 4 hits the capacity limits after the second step, so active power of Generator 2 is increased (Fig. 7); after 7 steps, Generator 2 hits the capacity limits. Since then, the re-dispatch capacity only comes from the dispatch of $\Delta\boldsymbol{x}[t]_{\mathrm{O}}$, which is very small and can be ignored. Therefore, during the second load peak, there is no remaining re-dispatch capacity for damping improvement, which is consistent with the results in Fig. 6.

In practice, there may be 20% extra re-dispatch capacity from the spinning reverse that is not considered in the above simulation. With this extra re-dispatch capacity, damping of the inter-area mode can be further improved, which is validated by the result in Fig. 6.

In the final part, we validate the damping improvement by setting a three-phase ground fault at bus 8 at time 11:00am. The relative angle of the Generators 1 and 3 is centered, normalized and shown in Fig. 9 for an illustrative comparison. Clearly, with the proposed method, damping of the inter-area mode is increased significantly; with extra re-dispatch capacity, the damping is further improved; the results are consistent with those in Fig. 6.

To sum up, all the results in Figs. 6 to 9 validate the effectiveness of the proposed method.

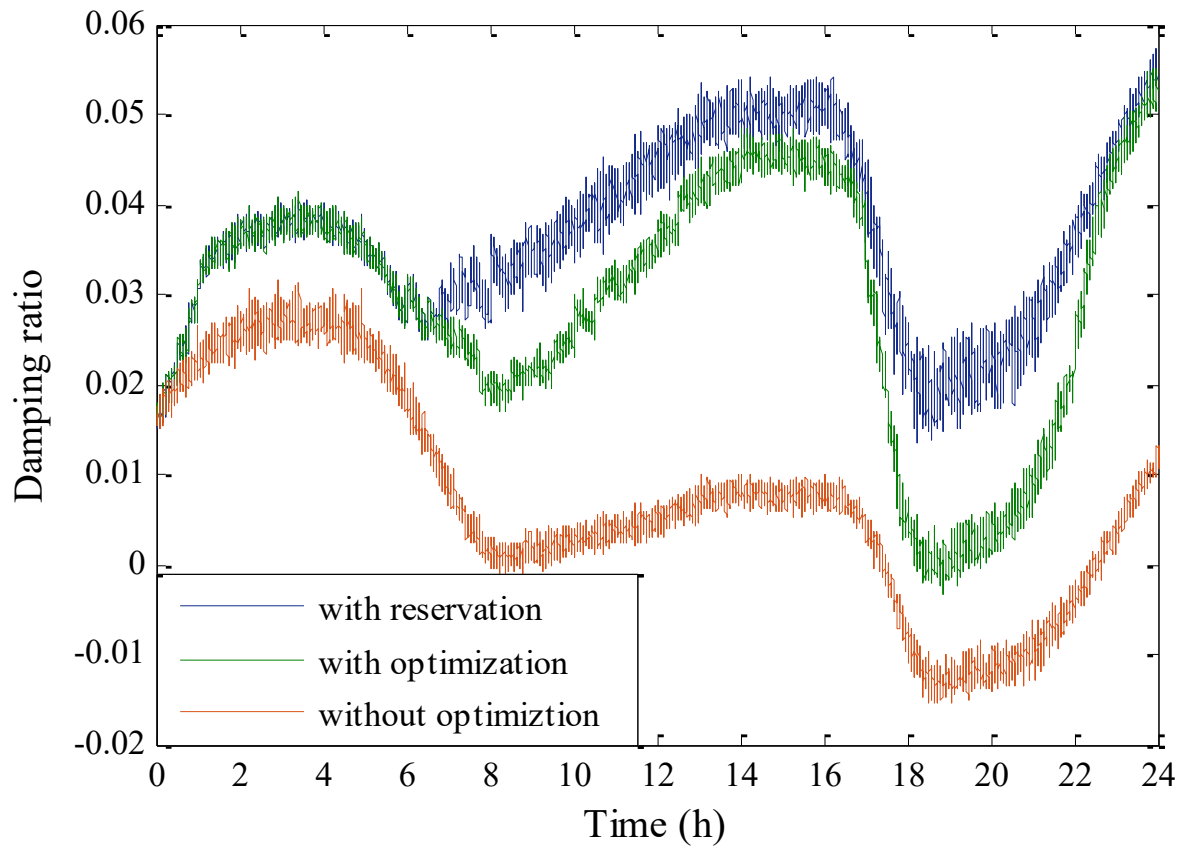

Fig. 6. Damping of the inter-area mode in three different manipulations.

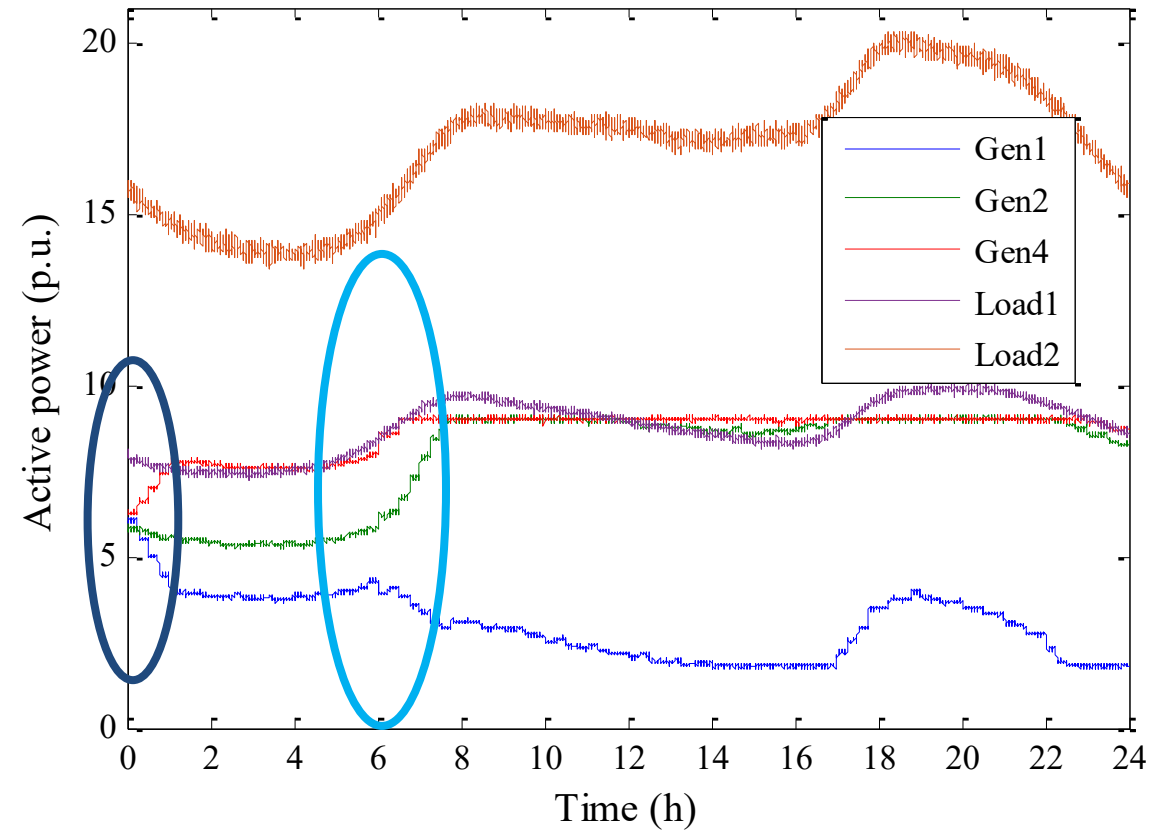

Fig.7. Active power flow of the generators and loads after optimization

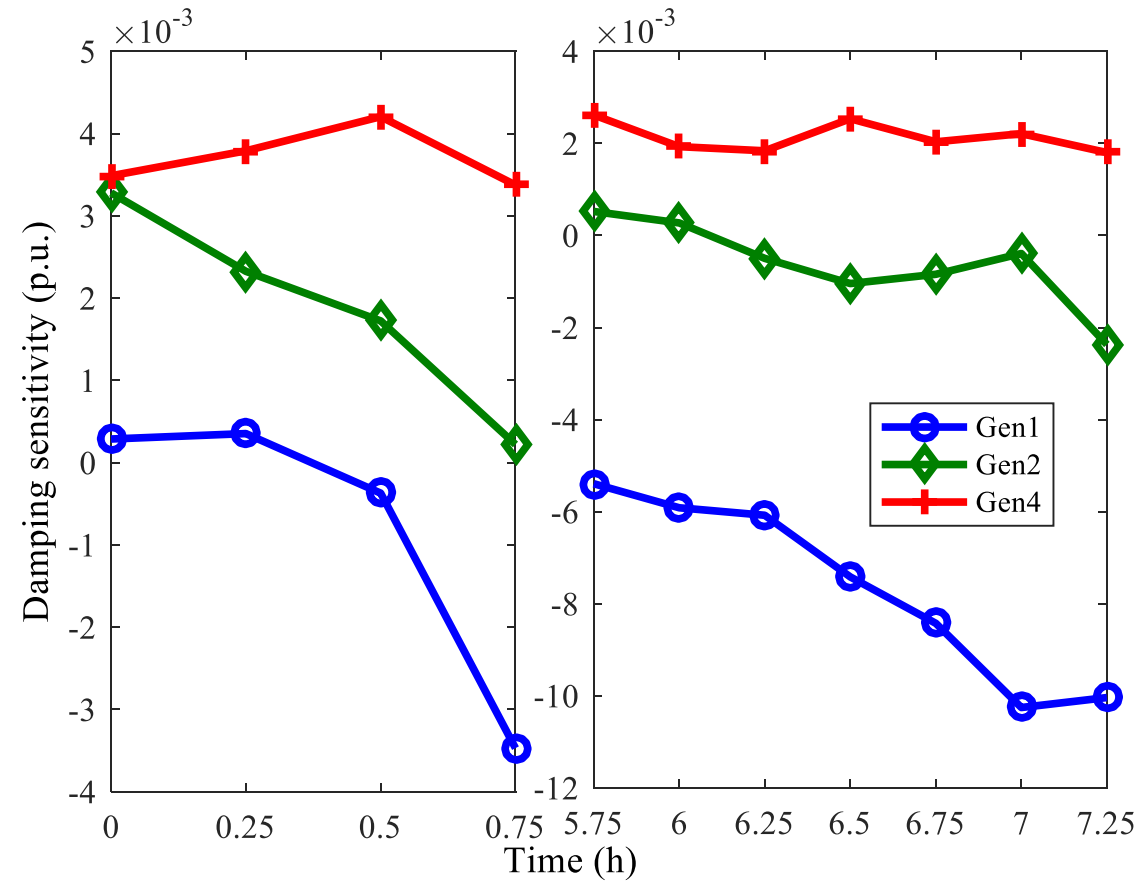

Fig. 8. Damping sensitivity during the sequential optimization.

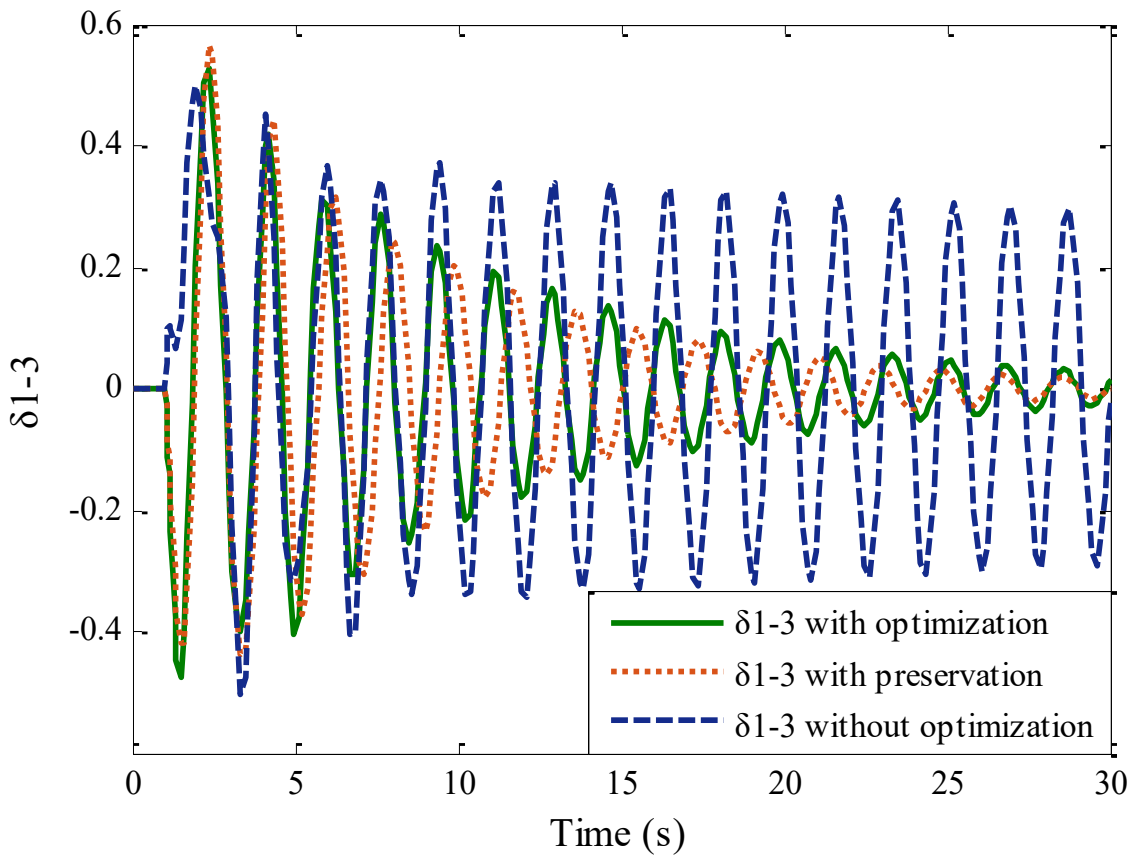

Fig. 9. The relative angle of Generators 1 and 3 in the 11 o'clock fault.

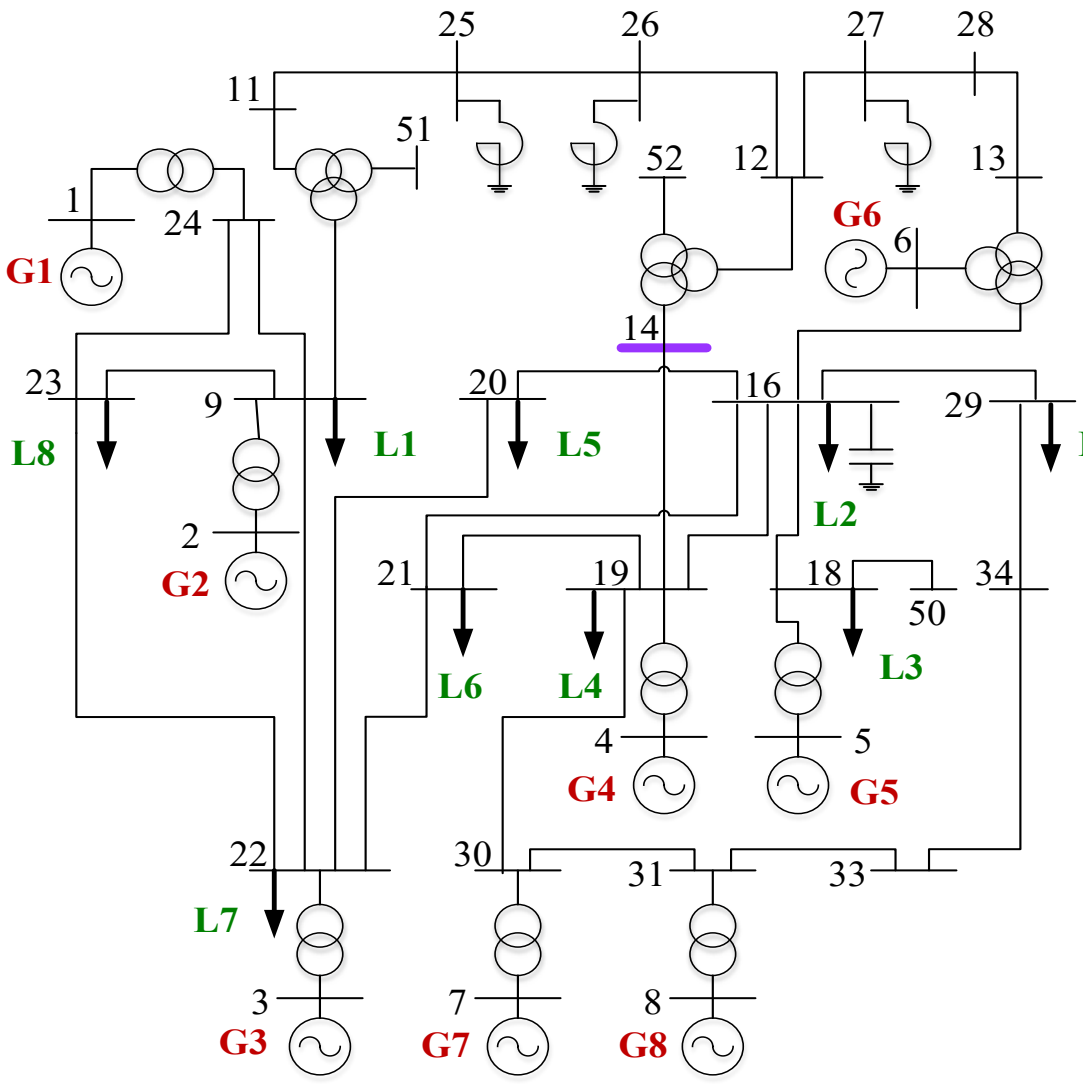

Fig. 10. The 8-generator 36-node system.

## B. 8-Machine 36-Bus System

The 8-generator 36-node system in Fig. 10 is used in this subsection. PSAT is used for simulation and the system condition is the same as the one in [30], to which readers are referred for details. In the power flow calculation, Generators 2-5, 7, and 8 are considered as P-V buses, with their active power as the controllable variables. Generator 1 is the slack bus and Generator 6 is a synchronous condenser. All loads are considered as P-Q buses. The original system has 7 electro-mechanical modes and the dominant mode is a poorly damped inter-area mode which dominates the oscillation between Generator 1 and Generators 3, 5, 7, and 8 [25]. In the following part, we validate the capability of the proposed approach to improve the damping of the dominant inter-area mode.

The three parallel threads in Section II Part C are employed.

1) Parallel thread 1:

The same method applied to the 2-area 4-machine system has been employed. The difference is that the fluctuated data is interpolated to 1 second per sample to create more samples during a dispatch interval. And $\underline{\boldsymbol{x}}[t+1]$ is set to 5% of the generator capacity here. The final generation curves ($\boldsymbol{x}[t]$) are shown in the upper plot of Fig. 11.

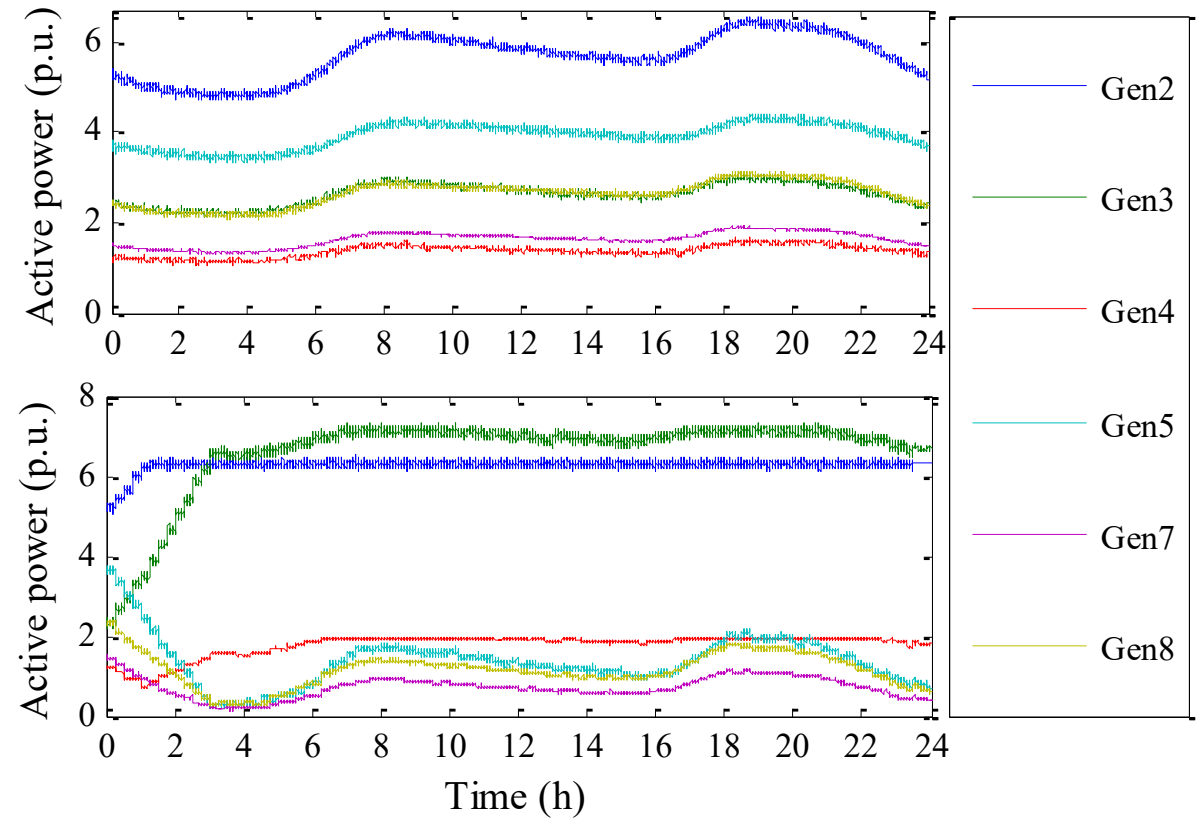

Fig. 11. Active power flow of the generators before and after optimization.

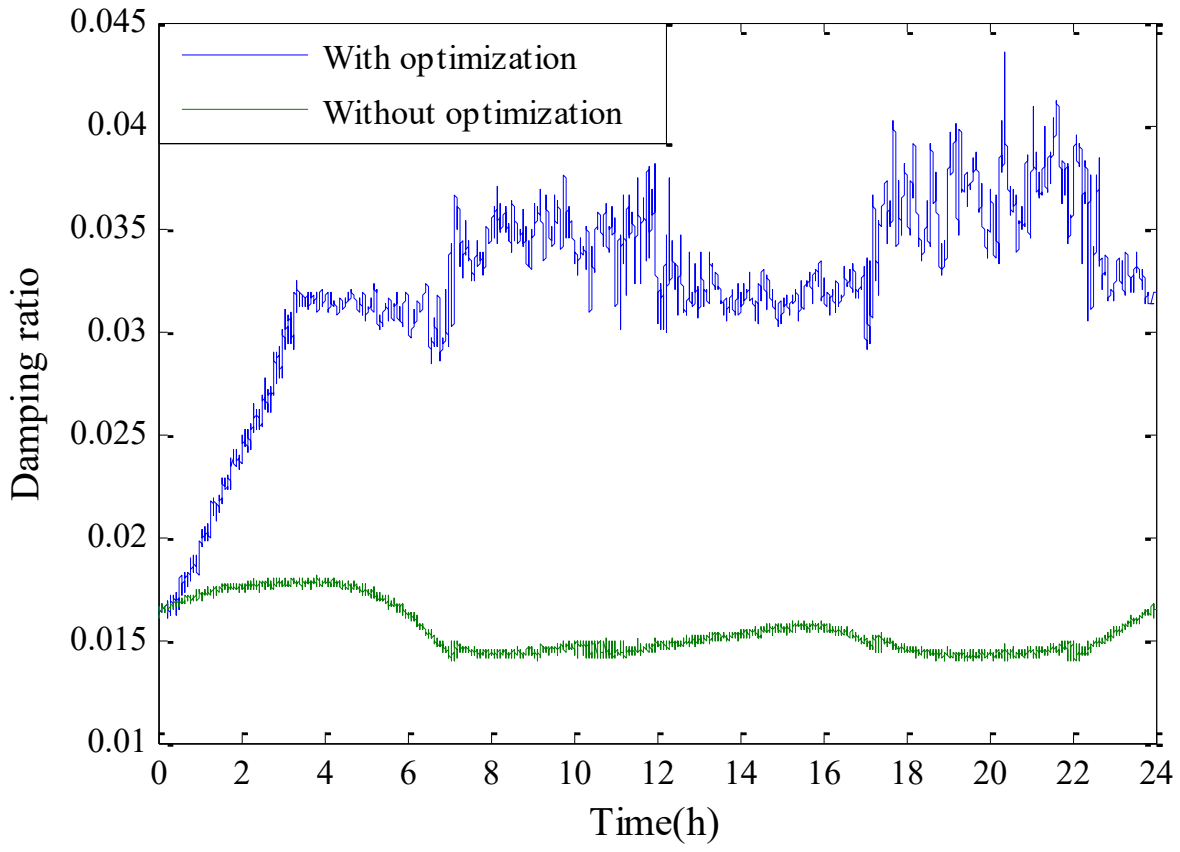

Fig. 12. Damping of the electromechanical mode 1.

1) Parallel thread 2:

Damping of the dominant inter-area mode is calculated using the same method applied to the 2-area 4-machine system. Its sampling rate is also enhanced to 1Hz, which could be achieved with recursive methods [13]. The damping curve ($\zeta[t]$) varying with the system power flow is shown in Fig. 12. Although each load and generator has limited impact on the dominant inter-area mode, without optimization, its damping is always under 3%.

2) Parallel thread 3:

The initial damping of the dominant inter-area mode is 1.64%; therefore, according to the trigger condition in parallel thread 3, the online optimization is carried out at the subsequent point of time. The constraints are the same as that in the 2-area 4-machine system. In the sensitivity estimation, $k$ = 0, the number of ensemble equals to 100. Table II shows $\boldsymbol{\Psi}[t]$ at the first re-dispatch point; the estimated sensitivity is similar to the one calculated by the model-based perturbation method which validates the effectiveness of the noise assisted ensemble regression.

Once the optimization is triggered, (5) is sequentially solved. The active generation after the re-dispatch is shown in the lower plot of Fig. 11, and the improved damping curve is compared with the initial one in Fig. 12. Clearly, after the first round of sequential optimization, damping of the domi-nant inter-area mode is improved beyond 3%. During the load peak, damping of the dominant inter-area mode tends to reduce (see the damping curve without optimization), while with several sequential re-dispatches, it keeps beyond 3%.

The damping sensitivity during the first round of sequential optimization is shown in Fig. 13. It can be figured out, the sensitivity results are consistent with the optimization based re-dispatch; although some of the generators hit the limits during the sequential optimization, after 12 steps, damping of the inter-area mode is improved over 3% (Fig. 12) and then the re-dispatch stops temporarily.

Moreover, we validate the damping improvement by setting a three-phase ground fault at bus 16 at time 20:00. The relative angle of the Generators 8 and 1 is centered, normalized and shown in Fig. 14. Clearly, with the proposed method, damping of the inter-area mode is increased significantly.

TABLE II
$\Psi$ AT THE FIRST RE-DISPATCH POINT

| Variable no. | 1 | 2 | 3 | 4 | 5 |
|---|---|---|---|---|---|
| Model-based | 0.0016 | 0.0011 | 0.0009 | -3.88E-5 | 0.0004 |
| NAER | 0.0018 | 0.0012 | 0.0010 | -2.81E-5 | 0.0005 |
| Variable no. | 6 | 7 | 8 | 9 | 10 |
| Model-based | 0.0003 | -0.0013 | -0.0013 | -0.0011 | -0.0013 |
| NAER | 0.0004 | -0.0014 | -0.0014 | -0.0012 | -0.0012 |
| Variable no. | 11 | 12 | 13 | 14 | 15 |
| Model-based | -0.0013 | -0.0011 | -0.0015 | -0.0014 | -0.0012 |
| NAER | -0.0013 | -0.0013 | -0.0013 | -0.0015 | -0.0014 |

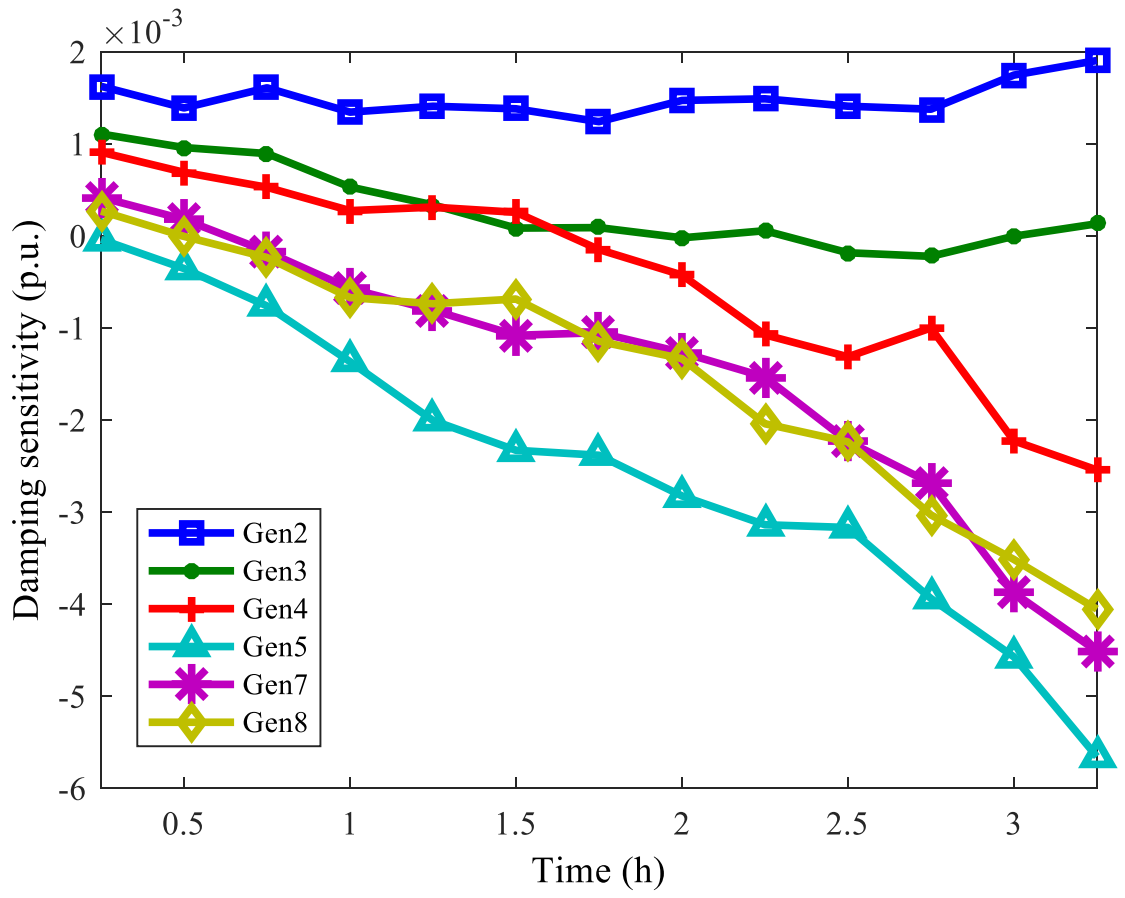

Fig. 13. Damping sensitivity during the sequential optimization.

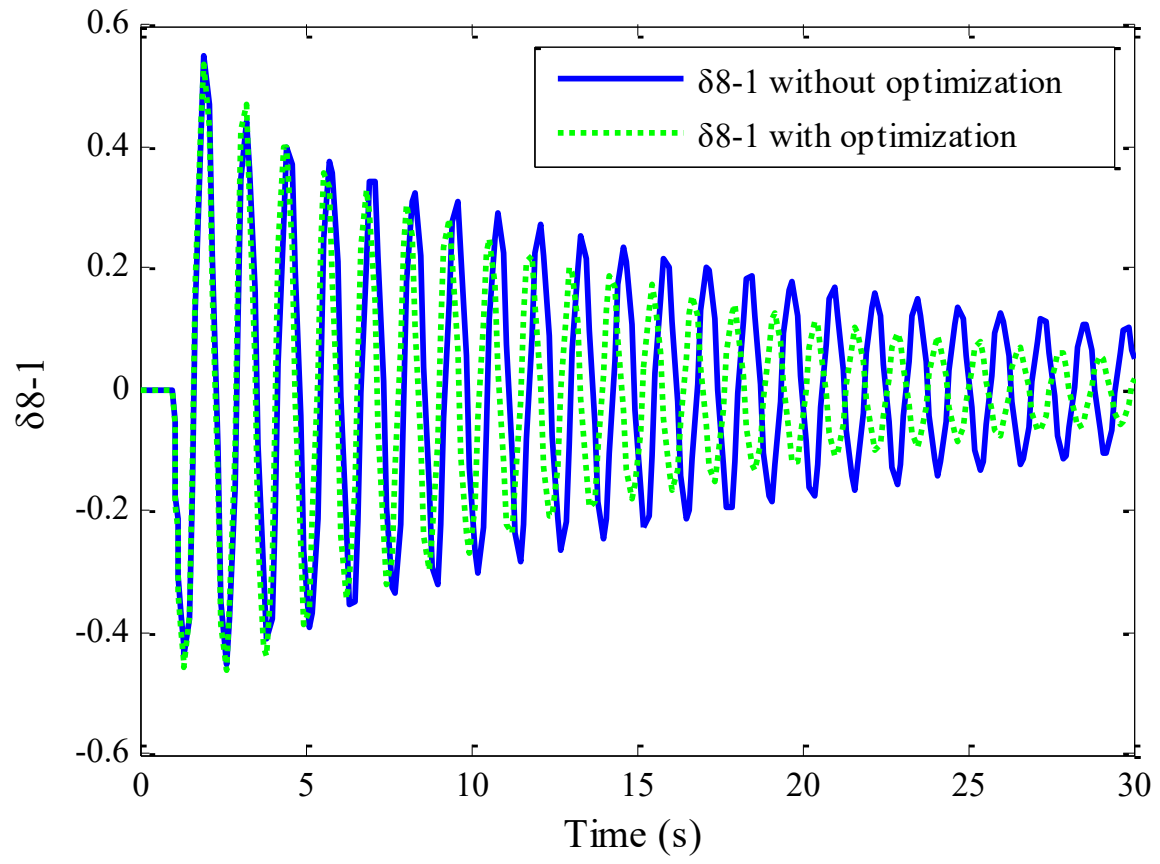

Fig. 14. The relative angle of Generators 8 and 1 in the 20 o'clock fault.

To sum up, the results in Figs. 11 to 14 validate the effectiveness of the proposed method in a relative larger system where each generator has limited impacts to the dominant inter-area mode.

## C. *NETS-NYPS 68-Bus System*

In this subsection, the proposed approach is validated by using the NETS-NYPS 68-bus system (Fig. 15), modeled from the realistic system – New England /New York power system, which has been widely used for small signal stability analysis and control and has been reported in many existing works [27]. Its static and dynamic simulations are performed in PSAT. In the power flow calculation, Generators 1-12 and 14-16 are considered as P-V buses, with their active power as the controllable variables, while generator 13 is the slack bus. All loads are considered as P-Q buses. Each generator is equipped with exciter and only Generator 9-13 are equipped with PSSs. The original system has several electromechanical modes and the dominant mode is a poorly damped inter-area mode dominates the oscillation between Generator 13 and the other generators. In the following part, we validate the proposed approach to improve damping of the dominant inter-area mode. Due to the limitation of paper space, we just put on the final result.

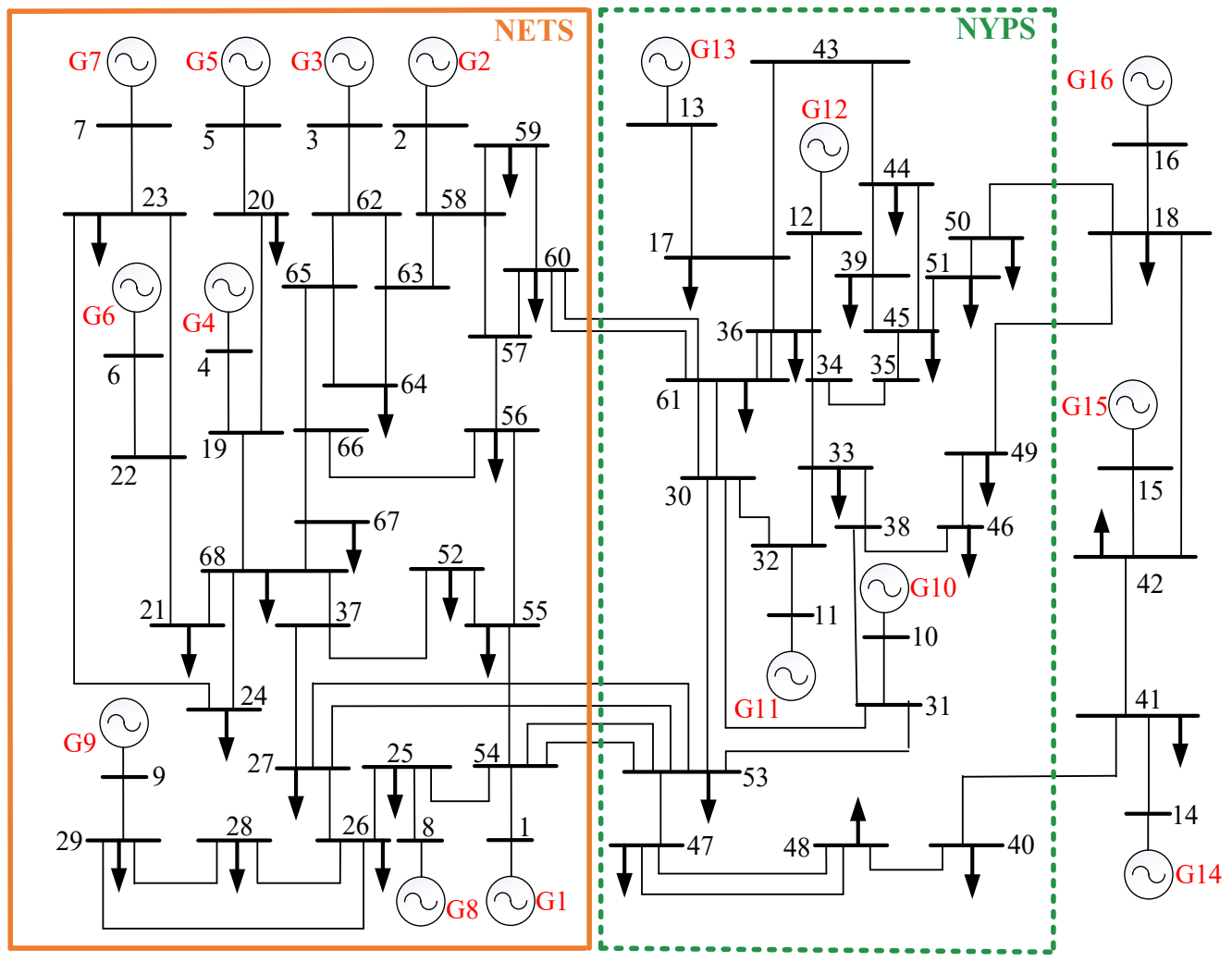

Fig.15. NETS-NYPS 68-bus system.

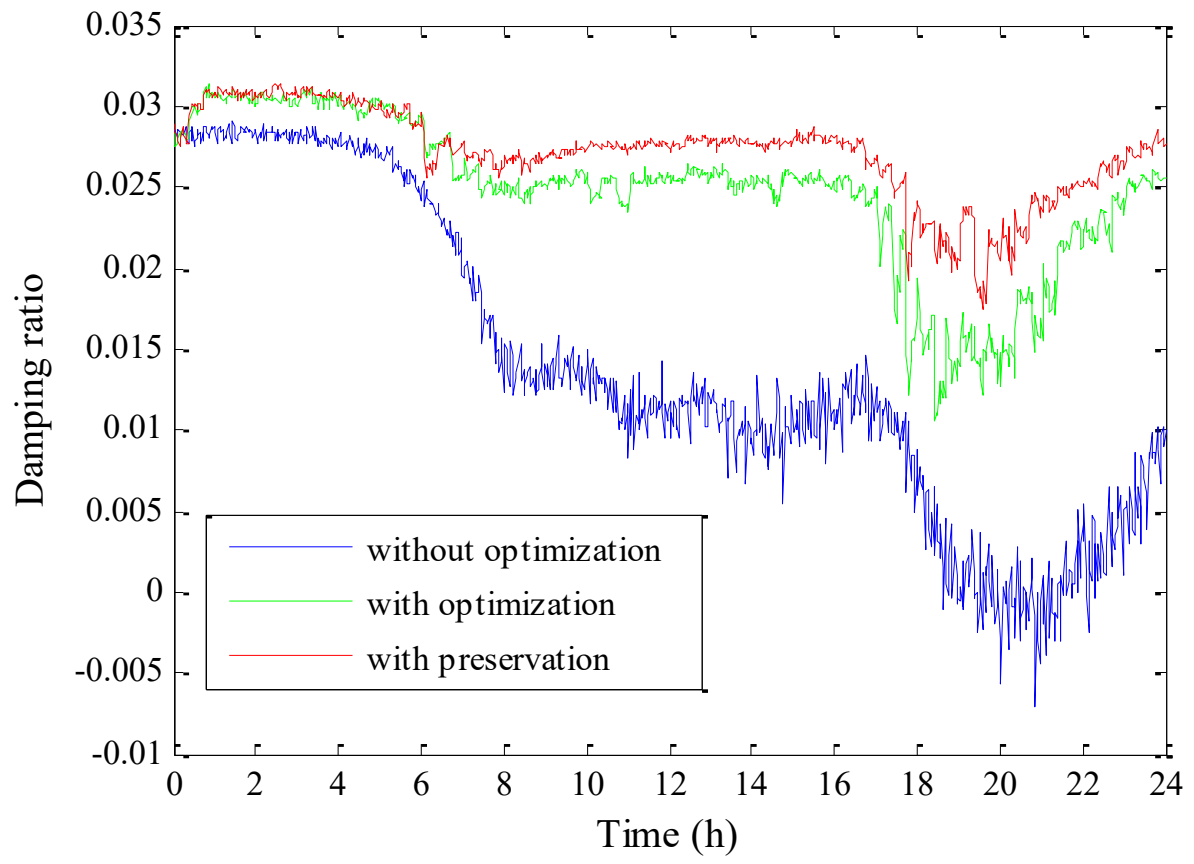

Fig. 16. Damping of the electromechanical mode.

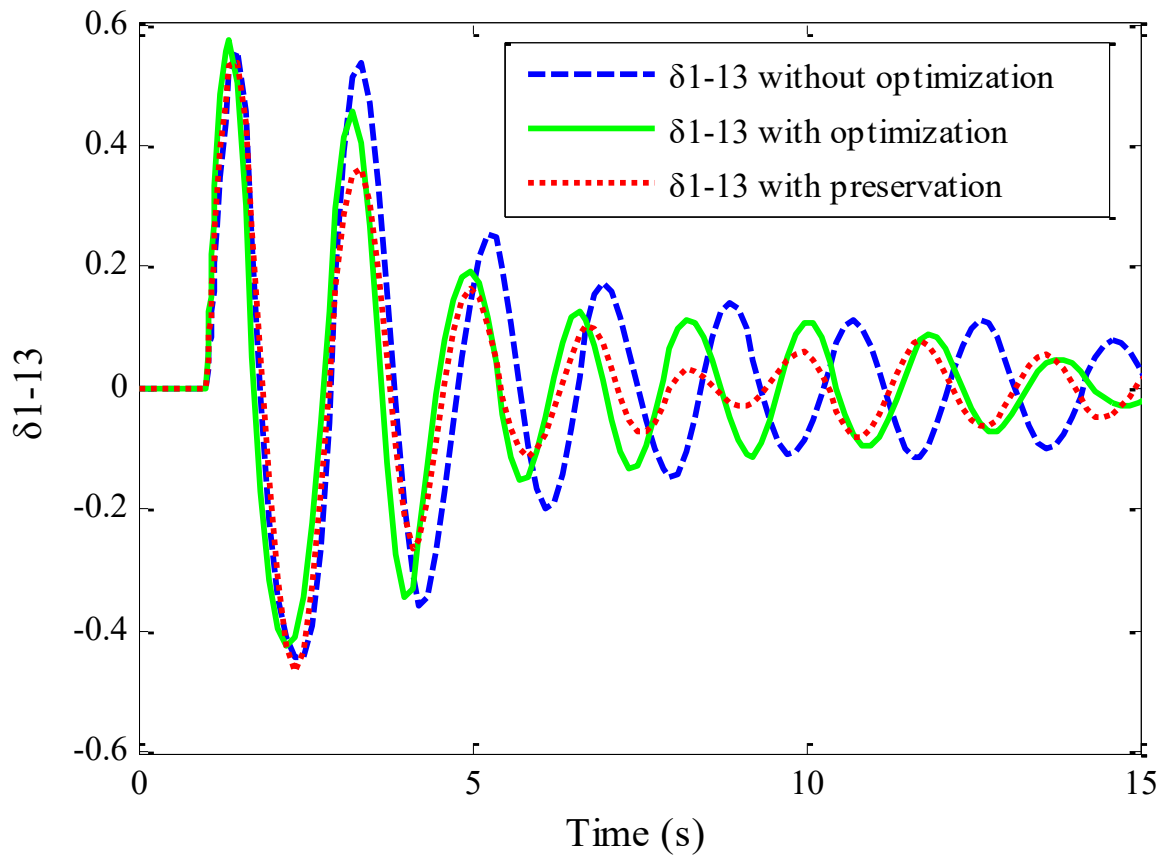

Fig. 17. The relative angle of Generators 1 and 13 in the 20 o'clock fault.

The three parallel threads in Section II Part C are employed. It is quite the same as the cases before. The great difference is that the interval time for adjustment here is 20 minutes since the increase number of related generators.

Once the optimization is triggered, (5) is sequentially solved. With the re-dispatch strategy, damping of the inter-area mode is improved to 3% before 6:00 and the system is stable during the second load peak. With this extra re-dispatch capacity, damping of the inter-area mode can be further improved and is validated by the result in Fig. 16. A three-phase ground fault is set at bus 60 at time 20:00. The relative angle of the Generators 1 and 13 is centered, normalized and shown in Fig. 17. Clearly, with the proposed method, damping of the inter-area mode is increased, which validated the effectiveness of the proposed method in a real system.

## V. Conclusion

This paper reports an initial work on power system oscillation damping improvement using data-driven online optimization method. Three simulation cases illustrate that the proposed optimization method is triggered once the damping of the system dominant mode is estimated to fall below 3%, and then an optimized re-dispatch that is compatible with the dispatch plans is sequentially carried out to improve the damping beyond 3%. The optimization constraints during the sequential re-dispatch are calculated with online sensitivity estimation and system state measurements; therefore, the proposed method is purely data-driven and can track the system operation condition in real-time without the need for the system offline models. The validated advantages show great potentials of the proposed method for application in real bulk power systems.

**Zhihao Chen** (S'17) was born in Guangdong, China, on September 25, 1995. He received his B.Eng. degree in electrical engineering from the School of Electric Power, South China University of Technology, in 2018. He has studied power system analysis, statistics, and several power and energy engineering related courses. His research interests include power system stability and control, system identification and data analysis in power system.

**Hanchen Xu** (S'12) received the B.Eng. and M.S. degrees in electrical engineering from Tsinghua University, Beijing, China, in 2012 and 2014, respectively, and the M.S. degree in applied mathematics from the University of Illinois at Urbana-Champaign, Urbana, IL, USA, in 2017, where he is currently working toward the Ph.D. degree at the Department of Electrical and Computer Engineering. His current research interests include optimization, reinforcement learning, with applications to power systems and electricity markets.

**Junbo Zhang** (S'10-M'14-SM'19) received his B.Eng. and Ph.D. from Tsinghua University in 2008 and 2013, respectively, followed by a postdoc fellow at the same university. He studied at The Hong Kong Polytechnic University from 2009 to 2010, and is now visiting Stanford University for a one year program. He is a professor at the School of Electric Power, South China University of Technology. His research areas include power system stability and control, knowledge management and intelligent decisions in power systems and power system economics.

**Lin Guan** (M'14) received the B.Eng. degree and the Ph.D. degree in electrical engineering from Huazhong University of Science and Technology, China, in 1990 and 1995, respectively. She is now a professor and the chair in the School of Electric Power, South China University of Technology. Her research interests include power system stability and control, planning and operation, and smart grid